\documentclass[12pt, reqno]{amsart}
\usepackage{amsmath}
\usepackage{amssymb}
\input amssym.def
\input amssym.tex

\font\tengothic=eufm10
\font\sevengothic=eufm7
\newfam\gothicfam
      \textfont\gothicfam=\tengothic
      \scriptfont\gothicfam=\sevengothic

\numberwithin{equation}{section}
\newenvironment{mydef}{{\bf Definition.}}{\hspace*{\fill} \par\vspace{1ex}}

\begin{document}

\setlength{\baselineskip}{1.5em}
\setcounter{section}{-1}
\newtheorem{thm}{\bf Theorem}

\newcommand{\pp}{{\mathbb P}}
\newcommand{\px}{\pp^2(\x)}
\newcommand{\zz}{{\mathbb Z}}
\newcommand{\nn}{{\mathbb N}}
\newcommand{\x}{{\mathbb X}}
\newcommand{\y}{{\mathbb Y}}
\newcommand{\ix}{I_{\mathbb X}}
\newcommand{\cv}{{\mathcal V}}
\newcommand{\calo}{{\mathcal O}}
\newcommand{\cl}{{\mathcal L}}
\newcommand{\cf}{{\mathcal F}}
\newcommand{\ck}{{\mathcal K}}
\newcommand{\ci}{{\mathcal I}}
\newcommand{\sfrac}[2]{\frac{\displaystyle #1}{\displaystyle #2}}

\title{Asymptotic $N_p$ property of rational surfaces}
\author{H\`a Huy T\`ai}
\address{Institute of Mathematics, P.O. Box 631, B\`o H\^o, H\`a N\^oi 10000, Vietnam} 
\email{tai@hanimath.ac.vn}
\thanks{{\it 2000 Mathematics Subject Classification.} 14E25, 14J26, 13D02}
\keywords{Rational surfaces, blowing up, fat points, resolution, $N_p$ property}
\maketitle

\section{Introduction.}

The pioneering work of Mumford (\cite{m}), its amplifications by St. Donat (\cite{sd}) and Fujita (\cite{f}), the inspiring work of Green (\cite{gr}), followed by works of Green and Lasarsfeld (\cite{grla}) and Ein and Lasarsfeld (\cite{el}), have captured the interest and have influenced a large number of researchers in the last fifteen years. Several authors have studied the defining equations of projective varieties and, more generally, the higher order syzygies among these equations. A significant algebraic property was introduced along this line of works (\cite{gr}, \cite{grla}), the {\it $N_p$ property}. It says that a variety is generated by quadratics, and its minimal free resolution is linear up to the first $p$ steps. We recall in detail the definition of $N_p$ property from \cite{grla}. 

\begin{mydef}
Let $Y$ be a smooth projective variety and let $\cl$ be a very ample line bundle on $Y$ defining an embedding $\varphi_{\cl}: Y \hookrightarrow \pp = \pp(H^0(Y, \cl)^{*})$. Let $S = \mbox{Sym}^{*} H^0(Y, \cl)$, the homogenenous coordinate ring of the projective space $\pp$. Suppose $A$ is the homogeneous coordinate ring of $\varphi_{\cl}(Y)$ in $\pp(H^0(Y, \cl)^{*})$, and 
\[ 0 \rightarrow F_n \rightarrow F_{n-1} \rightarrow \ldots \rightarrow F_0 \rightarrow A \rightarrow 0 \]
is a minimal free resolution of A. The line bundle $\cl$, or the embedding $\varphi_{\cl}(Y)$ of $Y$, is said to have {\it property $N_p$} (for $p \in \nn$) if and only if $F_0 = S$ and $F_i = S(-i-1)^{\alpha_i}$ with $\alpha_i \in \nn$ for all $1 \le i \le p$. One also says that $\cl$, or the embedding $\varphi_{\cl}(Y)$ of $Y$, satisfies property $N_0$ if $\varphi_{\cl}(Y)$ is projectively normal, i.e. $\cl$ is normally generated. 
\end{mydef}

The question of having property $N_p$, for various values of $p$, has been studied on many different projective varieties, such as the Veronese varieties (cf. \cite{op}, \cite{jpw}), and the Segre product of several projective spaces (cf. \cite{gp}, \cite{jpw}, \cite{r}). In this paper, we look at the same question for rational surfaces obtained by blowing up $\pp^2$ at collections of points. Defining equations and the minimal free resolution of these surfaces have also been the objects of study for many authors, such as \cite{cs}, \cite{chtv}, \cite{ch}, \cite{gg}, \cite{ggh}, \cite{ggp}, \cite{gi1}, \cite{gi2}, \cite{gl}, \cite{h}, \cite{hir}, \cite{ho1}, \cite{ho2}, \cite{ha1}, \cite{ha2}. Our work in this paper owes much of its origin to the work of \cite{gr} and \cite{ggp}.

Our main result is the following theorem.

{\bf Theorem.} {\it Suppose $\x = \{ P_1, \ldots, P_s \}$ is a set of $s$ distinct points in $\pp^2$. $Z = m_1P_1 + \ldots + m_sP_s \subseteq \pp^2$ is an arbitrary scheme of fat points with support $\x$ ($m_i \in \nn \ \forall i$). Let $\px$ be the blowup of $\pp^2$ along the points in $\x$, $D_t = tE_0 - \sum_{i=1}^s m_iE_i \in \mbox{Pic }(\px)$, and $\sigma = \sigma(Z)$. Then, for each $p \in \zz_{\ge 0}$, $\calo_{\px}(D_t)$ has property $N_p$ for all $t \ge \max \{\sigma+1, d, 1+\sfrac{d+p}{3}\}$, where $d = \sum_{i=1}^s m_i$. In other words, the embedding of $\px$ using $D_t$ satisfies property $N_p$ for all $t \ge \max\{\sigma+1, d, 1+\sfrac{d+p}{3}\}$. }

This theorem gives a positive answer to a conjecture of Geramita and Gimigliano (\cite{gg}), where they considered the case when $Z$ is reduced (i.e. $m_i = 1$ for all $i=1, \ldots, s$) and stated that the embedding of $\px$ using $D_t$ should have property $N_p$ for $t$ big enough ($t \ge \sigma(Z) + p$). In this situation, even though our bound is not exactly the conjectural value $\sigma(Z) + p$, it will be easily seen that for $p \gg 0$, $\max\{\sigma+1, d, 1+\sfrac{d+p}{3}\} \ll \sigma+p$. 

\section{Preliminaries and Notations.}

Suppose $\x = \{ P_1, \ldots, P_s \} \subseteq \pp^2 = \pp^2_k$ is a set of $s$ distinct points ($k$ is an algebraic closed field of characteristic 0) and let $\wp_i \subseteq R = k[w_0, w_1, w_2]$ be the defining ideal of $P_i$ for all $i = 1, \dots, s$. We also let $\px$ be the blowup of $\pp^2$ along the points in $\x$. It is well-known (see \cite{hart}) that $\mbox{Pic }(\px) \simeq \zz^{s+1} \simeq \ <E_0, E_1, \ldots, E_s>$, where $E_1, \ldots, E_s$ are the classes of the exceptional divisors and $E_0$ denotes the class of the pull-back of a general line in $\pp^2$. 

An approach to the study of the embeddings of $\px$ that has been investigated is through ``fat points" in $\pp^2$ with support $\x$. To be more precise, suppose $m_1, \ldots, m_s$ are positive integers, and let $I = \wp_1^{m_1} \cap \ldots \cap \wp_s^{m_s} \subseteq R$. Then, the {\it scheme of fat points} associated to $I$ is the subscheme $Z = m_1P_1 + \ldots + m_sP_s$ of $\pp^2$ defined by the ideal $I$. Let $\sigma = \sigma(I) = \sigma(Z)$ be the least integer start from which the difference function of the Hilbert function of $I$ vanishes (cf. \cite{cgo}). For each $t$, we have a divisor $D_t = tE_0 - \sum_{i=1}^s m_iE_i$ on $\px$. By a general result of Coppens (cf. \cite{ggp}) or a stronger version in $\pp^2$ of Davis and Geramita (\cite{dg}), $D_t$ is very ample on $\px$ for all $t \ge \sigma+1$, and $D_{\sigma}$ is also very ample provided any line in $\pp^2$ meets $Z$ at less than $\sigma$ points (counted properly). It was shown in \cite{ggp} that the embedding of $\px$ using $D_t$, for $t \ge \sigma$ and $D_t$ is very ample, is projectively normal and arithmetic Cohen-Macaulay (a.CM), i.e. the embedding of $\px$ using $D_t$ satisfies property $N_0$. When the scheme $Z$ is reduced (i.e. $m_i = 1$ for all $i$) it was shown in \cite{gi2} and \cite{gg} that the embedding of $\px$ using $D_t$, for $t \ge \sigma+1$, not only is projectively normal and a.CM, but also is generated by quadratics (i.e. satisfying properties $N_0$ and $N_1$). In \cite{gg}, the authors further believed and conjectured that for $t \ge \sigma+p$, the embedding of $\px$ using $D_t$ should have property $N_p$.

We shall now recall some notation and result introduced by Green (\cite{gr}) that will be used for our proof of the main theorem. 

Let $Y$ be a projective scheme. Let $\cl$ be a very ample line bundle and $\cf$ a coherent sheaf on $Y$. Let $W = H^0(Y, \cl)$ and $S = \mbox{Sym}^{*}W$. Then, $S$ is the homogeneous coordinate ring of $\pp(W)$, the projective space into which $Y$ is embedded using $\cl$. Let $B = B(\cl, \cf) = \bigoplus_{q \in \zz} H^0(Y, \cf \otimes q\cl) = \bigoplus_{q \in \zz} B_q$ a $S$-graded module. 

\begin{mydef} The {\it Koszul complex} of $B$ is defined to be
\[ \ldots \rightarrow \wedge^{p+1}W \otimes B_{q-1} \stackrel{d_{p+1, q-1}}{\longrightarrow} \wedge^pW \otimes B_q \stackrel{d_{p,q}}{\longrightarrow} \wedge^{p-1}W \otimes B_{q+1} \rightarrow \ldots \]
and the {\it Koszul cohomology groups} of $B$ are defined to be
\[ \ck_{p,q}(\cl, \cf) = \sfrac{\mbox{ker } d_{p,q}}{\mbox{im } d_{p+1, q-1}}, \ p,q \in \zz. \]
\end{mydef}

The following theorem relates the Koszul cohomology groups of $B$ and its minimal free resolution over $S$.

{\bf Theorem A.} (Green's syzygy theorem - \cite[1.b.4]{gr}) {\it Suppose 
\[ \ldots \rightarrow \bigoplus_{q \ge q_1} M_{1,q} \otimes S(-q) \rightarrow \bigoplus_{q \ge q_0} M_{0,q} \otimes S(-q) \rightarrow B \rightarrow 0 \]
is a minimal free resolution of $B$ over $S$, then
\[ M_{p, p+q}(\cl, \cf) = \ck_{p,q}(\cl, \cf) \ \forall p,q. \]}
Here, we write $M_{p,p+q}$ and $\ck_{p,q}$ as functions of $\cl$ and $\cf$ because $B$ itself depends on $\cl$ and $\cf$.

Green's notations and results are applicable to our situation when the scheme $Y$ is the blowup $\px$ of $\pp^2$, the line bundle $\cl$ is the invertible sheaf $\cl(D_t)$ corresponding to the divisor $D_t$ on $\px$. In this case, we write $M_{p,p+q}(D_t, \cf)$ and $\ck_{p,q}(D_t, \cf)$ for $M_{p, p+q}(\cl, \cf)$ and $\ck_{p,q}(\cl, \cf)$. When the coherent sheaf $\cf$ is the structure sheaf $\calo_{\px}$ of $\px$, we write $\ck_{p,q}(D_t)$ for $\ck_{p,q}(D_t, \cf)$. 

\section{Proof of the main theorem.}

For $p = 0$, the result was already proved by \cite{ggp}. Suppose $p \ge 1$. Let 
\[ d_p = \max \{ \sigma+1, d, 1+\sfrac{d+p}{3} \}, \]
where $d = \sum_{i=1}^s m_i$. 

Let $t$ be an arbitrary integer such that $t \ge d_p$. Let $D_t = tE_0 - \sum_{i=1}^s m_iE_i$ and $\cl = \cl(D_t)$ the invertible sheaf corresponding to $D_t$. Let $W = H^0(\px, \cl)$ and $S = \mbox{Sym}^{*} W$. Since $t \ge d_p \ge \sigma+1$, $D_t$ is very ample on $\px$. Let $\cv$ be the embedding of $\px$ into $\pp(W)$ using the divisor $D_t$. We need to show that $\cv$ possesses property $N_p$.

Let $I_V \subseteq S$ be the defining ideal of $\cv$, and let $\ci_V$ be the ideal sheaf of $\cv$ in $\pp(W)$. Since $D_t$ is very ample, we have an exact sequence
\[ 0 \rightarrow I_V \rightarrow S \rightarrow \bigoplus_{q \in \zz} H^0(\px, qD_t) \rightarrow \bigoplus_{m \in \zz} H^1(\pp(W), \ci_V(m)) \rightarrow 0. \]
Moreover, it was proved in \cite{ggp} that $\cv$ is projectively normal, i.e. 
\[ \bigoplus_{m \in \zz} H^1(\pp(W), \ci_V(m)) = 0; \]
or in other words, 
\[ S/I_V \simeq \bigoplus_{q \in \zz} H^0(\px, qD_t). \]
Thus, the minimal free resolution of the homogeneous coordinate ring of $\cv$ is given by that of $\bigoplus_{q \in \zz} H^0(\px, qD_t)$. This, by Green's syzygy theorem (Theorem A), is given by the Koszul cohomology groups $\ck_{r,n}(D_t), \ r,n \in \zz$. More precisely, let $N = \dim W -1$, then the minimal free resolution
\[ 0 \rightarrow F_{N-2} \rightarrow \ldots \rightarrow F_1 \rightarrow F_0 = S \rightarrow S/I_V \rightarrow 0 \]
of the homogeneous coordinate ring of $\cv$ (since $\cv$ is projectively normal and a.CM as proved by \cite{ggp}, and the codimension of $\cv$ is 2, the length of its minimal free resolution must be $N-2$) is given by
\[ F_i = \bigoplus_{q \ge 1} \ck_{i,q}(D_t) \otimes S(-i-q) \mbox{ for } i =1, \ldots, N-2. \]

From \cite[Corollary 2.6]{ggp}, we know that the Castelnuovo-Mumford regularity of $\ci_V$ is at most 3, so 
\[ \ck_{i,q}(D_t) = 0 \mbox{ for all } i=1, \ldots, N-2, \mbox{ and } q \ge 3. \]
Observe further that if we can show $\ck_{p,2}(D_t) = 0$ for all $p \in \nn$ and $t \ge d_p$, then the theorem is proved. This is because, by then, since $d_1 \le d_2 \le \ldots \le d_p$, we can use induction to show that $\cv$ has property $N_i$ for all $i = 1, \ldots, p$; in particular, it has property $N_p$.

It follows from \cite{ggp} that $H^1(\px, mD_t) = 0$ for all $m \in \zz$. Let $K_{\px}$ be the canonical divisor on $\px$. By Green's Duality theorem (\cite[2.c.6]{gr}), we have
\[ \ck_{p, 2}(D_t)^{*} = \ck_{N-2-p, 1}(D_t, K_{\px}). \]
Moreover, by Green's Vanishing theorem (\cite[3.a.1]{gr}), we have
\[ \ck_{N-2-p, 1}(D_t, K_{\px}) = 0 \mbox{ when } h^0(\px, K_{\px}+D_t) \le N-2-p. \]
Therefore, it remains to check that 
\begin{eqnarray} 
h^0(\px, K_{\px}+D_t) & \le & N-2-p \mbox{ for all } t \ge d_p. \label{cond}
\end{eqnarray}

It is well-known that $K_{\px} = -3E_0 + \sum_{i=1}^s E_i$. Thus, 
\[ K_{\px} + D_t = (t-3)E_0 - \sum_{i=1}^s (m_i-1)E_i. \]
We have $D_t^2 = t^2 - \sum_{i=1}^s m_i^2 > t^2 - (\sum_{i=1}^s m_i)^2 \ge 0$ since $t \ge d_p \ge d = \sum_{i=1}^s m_i$, so
\[ (K_{\px}+D_t).D_t = K_{\px}.D_t + D_t^2 > K_{\px}.D_t. \]
Therefore, $H^2(K_{\px}+D_t) = 0$ (\cite[Lemma V.1.7]{hart}). 

It follows from Serre's Duality theorem that $H^1(K_{\px}+D_t) = H^1(-D_t)$. By Kodaira's theorem (cf. \cite[p. 248]{hart}), since $D_t$ is very ample, we have $H^1(-D_t) = 0$. Thus, $H^1(K_{\px}+D_t) = 0$.

Now, using the Riemann-Roch theorem, we obtain 
\begin{eqnarray*}
h^0(\px, K_{\px}+D_t) & = & \sfrac{1}{2}(K_{\px}+D_t).D_t + 1 \\
& = & \sfrac{1}{2}\big( t(t-3) - \sum_{i=1}^s m_i(m_i-1) \big) + 1 \\
& = & \sfrac{1}{2} t(t-3) - \sum_{i=1}^s {m_i \choose 2} + 1.
\end{eqnarray*}
Moreover, the Hilbert function of $Z$ is increasing and stablizes at $\deg(Z) = \sum_{i=1}^s {m_i+1 \choose 2}$.
Therefore, 
\[ N \ge {t+2 \choose 2} - \deg(Z) -1 = \sfrac{1}{2} (t+1)(t+2) - \sum_{i=1}^s {m_i+1 \choose 2} -1. \]
Hence, 
\begin{eqnarray*} 
N-2-p - h^0(\px, K_{\px}+D_t) & \ge & 3t - 3 - p - \sum_{i=1}^s {m_i+1 \choose 2} + \sum_{i=1}^s {m_i \choose 2} \\
& = & 3t - 3 - p - \sum_{i=1}^s m_i \\
& = & 3t - 3 - p - d \ge 0. \
\end{eqnarray*}
The inequality (\ref{cond}) is verified, and the theorem is proved. \qed

{\bf Remark:} When $p \gg 0$, $1+\sfrac{d+p}{3}$ is clearly the larger value compared to $\sigma+1$ and $d$. Thus, for $p \gg 0$, $\max\{\sigma+1, d, 1+\sfrac{d+p}{3}\} = 1+\sfrac{d+p}{3}$, which is essentially smaller than $\sigma+p$.

\begin{small}
{\it Acknowledgement:} {\sf I would like to thank Prof. A.V. Geramita for suggesting me this problem as well as many other related questions on rational surfaces.}
\end{small}

\end{document}